\documentclass[10pt,reqno]{amsart}
\usepackage{amsmath,amssymb}
\usepackage[french]{babel}
\usepackage[latin1]{inputenc}
\usepackage{epsfig}
\usepackage[all]{xy}

\newtheorem{teo}{Th\'eor\`eme}[section]
\newtheorem{cor}[teo]{Corollaire}
\newtheorem{pro}[teo]{Proposition}
\newtheorem{lem}[teo]{Lemme}

\theoremstyle{remark}
\newtheorem{rem}[teo]{Remarque}
\newtheorem{exe}[teo]{Exemple}
\newtheorem{exes}[teo]{Exemples}

\begin{document}
\newcommand{\ii}{\'{\i}}  
\renewcommand{\r}{\hspace*{1 cm}}
\newcommand{\codim}{\mbox{codim\,}}
 
\newcommand{\cre}{{\rm Cr}(\plan)}
\newcommand{\aut}{{\rm Aut}}
\newcommand{\jon}{{\rm Jon}}
\newcommand{\bir}{{\rm Bir}}
\newcommand{\hir}{{\mathbb F}}
\newcommand{\pic}{{\rm Pic}}
\newcommand{\rank}{\rm rg\,} 
\newcommand{\base}{\rm Base} 
\newcommand{\pgl}{{\rm PGL}}
 
\newcommand{\bl}{{\mathbb B}l} 
\newcommand{\mor}{\mbox{Mor}}
\newcommand{\ara}{\ar@} 

\renewcommand{\div}{{\rm div}} 
\newcommand{\supp}{{\rm Supp}}
\newcommand{\ord}{{\rm ord}}
\newcommand{\adj}{{\rm Adj}}
\newcommand{\fix}{{\rm Fix}}
\newcommand{\ie}{c'est-\`a-dire}
\newcommand{\Ie}{C'est-\`a-dire}
\newcommand{\proj}{{\mathbb P}}
\newcommand{\plan}{{\mathbb P}^2}
\newcommand{\complex}{{\mathbb C}}
\renewcommand{\natural}{{\mathbb N}}
\newcommand{\real}{{\mathbb R}}
\newcommand{\rational}{{\mathbb Q}}
\newcommand{\integer}{{\mathbb Z}}
\newcommand{\grass}{{\mathbb G}}
\newcommand{\N}{{\mathbb N}}
\newcommand{\R}{{\mathbb R}}
\newcommand{\Q}{{\mathbb Q}}
\newcommand{\Z}{{\mathbb Z}}
\newcommand{\C}{{\mathbb C}}
\newcommand{\G}{{\mathbb G}}
\newcommand{\poids}{\mbox{poids}}
\newcommand{\wtilde}{\widetilde}
\newcommand{\cohz}{{\rm H^0}}
\newcommand{\cohu}{{\rm H^1}}
\newcommand{\cohd}{{\rm H^2}}

\newcommand{\te}{\text}
\newcommand{\lto}{\longrightarrow}
\newcommand{\tor}{\xymatrix{\ar@{-->}[r]&}}
\newcommand{\mapstor}{\xymatrix{\ar@{|-->}[r]&}}
\newcommand{\cali}{{\mathcal I}}
\newcommand{\cala}{{\mathcal A}} 
\newcommand{\cale}{{\mathcal E}}
\newcommand{\calt}{{\mathcal T}}
\newcommand{\caln}{{\mathcal N}}  
\newcommand{\calf}{{\mathcal F}}
\newcommand{\calg}{{\mathcal G}} 
\newcommand{\calh}{{\mathcal H}}
\newcommand{\cald}{{\mathcal D}}
\newcommand{\calo}{{\mathcal O}}
\newcommand{\calw}{{\mathcal W}}
\newcommand{\calr}{{\mathcal R}}
\newcommand{\calp}{{\mathcal P}}
\newcommand{\calu}{{\mathcal U}}
\newcommand{\call}{{\mathcal L}}
\newcommand{\calm}{{\mathcal M}}
\hyphenation{In-ter-sec-tion  in-ter-sec-tion as-so-ci ral-ment 
trans-for-ma-tion bi-ra-tio-n-nel-le en-gen} 

\title[ - Sur un th\'eor\`eme de Castelnuovo]{Sur un th\'eor\`eme de Castelnuovo}

\address{J\'er\'emy Blanc\\ Laboratoire J.A. Dieudonn\'e (UMR 6621 du C.N.R.S.)\\
Universit\'e de Nice Sophia Antipolis\\
Parc Valrose\\ 
06108 Nice cedex 2\\ France}\email{blancj@math.unice.fr}
\address{Ivan 
Pan\\ Instituto de Matem\'atica\\UFRGS\\av. Bento Gon\c calves
9500\\91540-000 Porto Alegre, RS\\Brasil}\email{pan@mat.ufrgs.br}\address{Thierry Vust\\ Universit\'e de Gen\`eve\\ 
Section de math\'ematiques\\ 
2-4 rue du Li\`evre\\
CP 64, 1211 Gen\`eve 4\\ Switerzland}\email{Thierry.Vust@math.unige.ch}
\maketitle 
\begin{center}
J\'ER\'EMY BLANC, IVAN PAN\footnote{Partiellement soutenu par le CNPq-Brasil et la Section de Math\'ematiques de l'Universit\'e de Gen\`eve} et THIERRY VUST
\end{center}

\begin{abstract}
Nous poursuivons l'\'etude faite par G. Castelnuovo en 1892 au sujet du groupe des transformations birationnelles du plan complexe qui fixent points par points une courbe de genre $>1$; nous nous servons comme lui des syst\`emes lin\'eaires adjoints de la courbe fixe. 

Nous d\'emontrons que ces groupes sont ab\'eliens, et qu'ils sont soit finis, d'ordre 2 ou 3, soit conjugu\'es \`a un sous-groupe du groupe de de Jonqui\`eres. Nous montrons \'egalement que ces r\'esultats ne se g\'en\'eralisent pas aux courbes de genre $\leq 1$.

\noindent{\sc Mots-cl\'es.} transformations de Cremona, transformations birationnelles, courbes fixes, courbes de genre sup\'erieur, syst\`emes lin\'eaires adjoints, transformations de de Jonqui\`eres.

\vspace{0.3 cm}

\selectlanguage{francais}
\noindent{\sc Abstract.} 
{We continue the study of G. Castelnuovo on the group of birational transformation of the complex plane that fix each point of a curve of genus $>1$; we use adjoint linear system of the curve as Castelnuovo does.

We prove that these groups are abelian, and that these are either finite, of order 2 or 3, or conjuguate to a subgroup of the de Jonqui\`eres group. We show also that these results do not generalise to curves of genus $\leq 1$.}


\vspace{0.3 cm}

\end{abstract}

\section{Introduction}
On d\'esigne par $\plan=\proj^2_{\complex}$ le plan projectif sur le
corps
$\complex$ des
nombres complexes; une transformation de Cremona de $\plan$ est une application
birationnelle $\plan\tor\plan$ et on dit qu'une telle transformation
est de de Jonqui\`eres si elle pr\'eserve un pinceau de droites. L'ensemble des transformations de Cremona forme un groupe, appel\'e groupe de Cremona.

\vspace{0.3 cm}

L'origine de ce travail est l'envie de comprendre le tr\`es beau
th\'eor\`eme de G. Castelnuovo \cite{Cas}:

\vspace{0.1 cm}

\emph{Se una transformazione Cremoniana fra due piani sovrapposti
muta in s\`e stesso ciascun punto di una curva irreduttibile C di
genere superiore ad 1, la transformazione o \`e riduttibile al tipo
Jonqui\`eres, oppure \`e ciclica di $2^\circ$, $3^\circ$ o $4^\circ$
grado}\footnote{Traduction litt\'erale: \emph{Si une transformation Cremonienne entre
deux plans superpos\'es envoie sur eux-mêmes chaque point d'une courbe
irr\'eductible $C$ de genre sup\'erieur \`a 1, la transformation ou bien est
r\'eductible au type Jonqui\`eres, ou alors est cyclique du
$2^\circ$, $3^\circ$ ou $4^\circ$} degr\'e.}.

\vspace{0.1 cm}

Il affirme en particulier: une transformation de Cremona d'ordre
infini qui fixe points par points une courbe irr\'eductible de genre
(g\'eom\'etrique) $>1$ est conjugu\'ee \`a une transformation de de Jonqui\`eres. C'est ce
r\'esultat qui est peut-être le plus int\'eressant puisque on n'a que peu
de prise sur les transformations d'ordre infini.

\vspace{0.5 cm}

Dans cette note nous d\'emontrons une version un peu plus pr\'ecise du th\'eor\`eme de Castelnuovo:


\begin{teo}
Soit $F:\plan\tor\plan$ une transformation de Cremona diff\'erente de
l'identit\'e  qui fixe points par points une courbe
irr\'eductible de genre $>1$. Alors $F$ est conjugu\'ee \`a une
transformation de de Jonqui\`eres ou bien $F$ est d'ordre 2 ou 3. De
plus, dans le premier cas, si $F$ est d'ordre fini, c'est une
involution.
\label{teo1}
\end{teo}

Rappelons les exemples suivants.

\begin{exes} 

(\cite{Hud}, \cite{God}, \cite{Coo}, \cite{SR}, \cite{BayBea},
\cite{Fer}, \cite{Bla})

a) Les involutions de Geiser fixent une courbe non hyperelliptique de
genre 3.

b) Les involutions de Bertini fixent une courbe non hyperelliptique
de genre 4 \`a mod\`ele lisse sur un c\^one quadratique.

c) Les involutions de de Jonqui\`eres fixent une courbe
hyperelliptique. 
\label{exes1}
\end{exes}

\begin{exe} 

(\cite{Fer}, \cite{Dol}, \cite{Bla}). Dans l'espace projectif \`a poids
$\proj(3,1,1,2)$ consid\'erons une surface lisse $S$ d'\'equation
$w^2=z^3+F_6(x,y)$ o\`u $F_6$ est homog\`ene de degr\'e 6: c'est un type
particulier de surfaces de Del Pezzo de degr\'e 1. La restriction de
$(w:x:y:z)\mapsto (w:x:y:\omega z)$, o\`u $\omega\neq 1$ est une racine
cubique de l'unit\'e, d\'efinit un automorphisme de $S$ d'ordre 3 dont
l'ensemble des points fixes contient une courbe irr\'eductible de genre
2.
\label{exes2}
\end{exe}

Voici encore un autre type d'exemple.

\begin{exe}
Soit $h\in\C[x]$ un polyn\^ome de degr\'e $2g+2$ sans racines multiples.
Notons $J_h$ le tore de $\pgl(2,\complex(x))$ image du sous-groupe 
\[T_h:=\Big\{\left(\begin{array}{ll} a_1&ha_2\\a_2&a_1\end{array}\right):
a_i\in\complex(x), a_1^2-h a_2^2\neq 0\Big\}\]
de ${\rm GL}(2,\complex(x))$.

Alors, pour tout $A\in T_h$, en d\'esignant par $a$ l'image de $A$ dans $J_h$, l'application rationnelle
$F_a:\complex^2\tor\C^2$
d\'efinie par
\[(x,y)\mapsto \left(x,\frac{a_1y+ha_2}{a_2y+a_1}\right)\]
est de de Jonqui\`eres et laisse fixe la courbe hyperelliptique
$C$ d'\'equation $(y^2=h(x))$. Lorsque $a_1=0$, on retrouve l'exemple
\ref{exes1}c).

On observe que $T_h$ est isomorphe au groupe multiplicatif $\C(C)^{*}$  du corps $\C(C)$ des fonctions rationnelles sur $C$ et donc que $J_h$ est isomorphe \`a $\C(C)^*/\C(x)^*$.
\label{exes3}
\end{exe}

\vspace{0.5 cm}

En fait, les exemples  \ref{exes1}, \ref{exes2} et \ref{exes3} d\'ecrivent tous les types de transformations \'etudi\'ees par
Cas\-tel\-nuo\-vo, et ceci m\^eme lorsqu'on \'etend la recherche \`a des sous-groupes du groupe de Cremona, comme l'\'enonce le th\'eor\`eme \ref{teo2} ci-dessous.

Si $p\in\plan$, on d\'esigne par $\jon_p(\plan)$ le groupe des
transformations de de Jonqui\`eres en $p$, \ie\ celles qui stabilisent le pinceau des droites par $p$. \'Evidemment, si $q\in\plan$, le groupe $\jon_q(\plan)$
est conjugu\'e \`a $\jon_p(\plan)$ dans le groupe de Cremona; pour all\'eger on \'ecrira $\jon(\plan)$.

\begin{teo}
Soit $M\neq \{id\}$ un sous-groupe du groupe de Cremona de $\plan$; supposons
qu'il existe une courbe irr\'eductible $C$ de genre $>1$ qui est
laiss\'ee fixe points par points par tous les \'el\'ements de $M$. Alors
$M$ est cyclique d'ordre 2 ou 3 engendr\'e par l'une des
transformations des exemples \ref{exes1}, \ref{exes2}, ou bien $M$
est conjugu\'e \`a un sous-groupe de $J_h$ (exemple \ref{exes3}); en particulier si $M$ est
infini, $C$ est hyperelliptique et $M$ est ab\'elien conjugu\'e \`a un sous-groupe
de $\jon(\plan)$.
\label{teo2}
\end{teo}

\vspace{0.3 cm}

Notre d\'emarche suit exactement celle de Castelnuovo qui consiste \`a montrer
qu'il existe un pinceau de courbes rationnelles ou elliptiques qui
est laiss\'e fixe par la transformation $F$ dans l'\'enonc\'e du th\'eor\`eme
\ref{teo1} (resp. par toute $\mu\in M$ dans \ref{teo2}). Ce pinceau est obtenu en
construisant les \emph{syst\`emes lin\'eaires adjoints successifs} de la
courbe de points fixes de $F$. Cette m\'ethode est \'egalement d\'ecrite dans les
livres de L. Godeaux \cite[chap. VIII, \S 2]{God} et de J. L.
Coolidge \cite[Book IV, chap. VII, \S 3, Thm. 14]{Coo}.

\bigskip

Lorsque $M$ est fini, l'existence de ce pinceau se d\'emontre aussi
avec des m\'ethodes plus modernes. En effet, pour commencer on observe
qu'il existe une surface rationnelle lisse $S$ et une application
birationnelle $\varphi:S\tor \plan$ telle que $\varphi^{-1}\mu\varphi$ est
un automorphisme bir\'egulier de $S$, ceci pour tout $\mu\in M$
(\cite{FerEin}, \cite{Dol}). On peut donc supposer que $M$ est un
sous-groupe d'automorphismes de $S$ et de plus qu'on
se trouve dans l'une des situations suivantes (\cite{Isk},
\cite{Dol}): 

1) il existe un fibr\'e en coniques $\pi:S\to\proj^1$ et un
homomorphisme 
\[M\to \aut(\proj^1), \  \mu\mapsto \overline{\mu}\] 
tel que $\pi\circ\mu=\overline{\mu}\circ\pi$;

2) $S$ est une surface de Del Pezzo. 

Dans le premier cas, puisque $M$ fixe une courbe $C$ de genre $\geq 1$,
on a $\overline{\mu}=id$ pour tout $\mu\in M$: autrement dit $M$ fixe
un pinceau de courbes rationnelles.

Dans le second cas, consid\'erons l'application rationnelle
anticanonique \begin{center}$\gamma:S\tor |-K_S|^{\vee}$\end{center} qui est \'equivariante.
Puisque aucun \'el\'ement du syst\`eme anticanonique ne peut contenir $C$
(voir lemme \ref{lem34}), la courbe $\gamma(C)$ engendre $|-K_S|^{\vee}$ et
par cons\'equent $M$ op\`ere trivialement au but, d'o\`u l'existence d'un
pinceau de courbes elliptiques laiss\'e fixe par $M$.

\bigskip

Par contre, si $M$ est infini, la m\'ethode des ``adjoints successifs"
reste d'actualit\'e.

\bigskip

Finalement, nous d\'emontrerons \`a la derni\`ere section (Proposition \ref{propellrat}) que les th\'eor\`emes \ref{teo1} et \ref{teo2} ne se g\'en\'eralisent pas aux courbes rationnelles ou elliptiques (de genre $0$ ou $1$).
\section{Syst\`eme adjoint}

On adoptera la convention suivante: toutes les surfaces consid\'er\'ees sont tacitement suppos\'ees projectives
lisses, rationnelles et connexes.

Soit $C$ une courbe irr\'eductible contenue dans la surface $X$ et
$\pi:Y\to X$ une r\'esolution plong\'ee des singularit\'es de $C$; notons
$\widetilde{C}$ la transform\'ee stricte  inverse de $C$ par $\pi$. Si
$h^0(\widetilde{C}+K_Y)>1$, autrement dit si le syst\`eme lin\'eaire
$|\widetilde{C}+K_Y|$ n'est pas vide ni r\'eduit \`a un seul diviseur,
on note $\adj(C)$ le syst\`eme lin\'eaire $\pi_*|\widetilde{C}+K_Y|$
priv\'e de ses \'eventuelles composantes fixes: c'est le \emph{syst\`eme
adjoint}
de $C$. Ce syst\`eme est ind\'ependant de la r\'esolution
choisie : soient en effet $\pi_i:Y_i\to X$ ($i=1,2$) deux r\'esolutions
de $C$; Consid\'erons alors un diagramme commutatif
\[\xymatrix{ &Y\ar@{->}[dl]_{\sigma_1}\ar@{->}[rd]^{\sigma_2}& \\
Y_1\ar@{->}[rd]_{\pi_1}& & Y_2\ar@{->}[ld]^{\pi_2}\\
 & X& }\]
o\`u $\sigma_i$ ($i=1,2$) sont des morphismes birationnels, et notons
$\widetilde{C_i}$, $\widetilde{C}$ les transform\'ees strictes inverses
de $C$ par $\pi_i$ et $\pi_i\circ\sigma_i$ respectivement. Soit
maintenant $D_i\in|\widetilde{C_i}+K_{Y_i}|$ et supposons pour
simplifier que $\sigma_i$ est l'\'eclatement d'un point $O_i$; alors le
transform\'e total $\sigma_i^* D_i$ de $D_i$ appartient \`a
$|\wtilde{C}+m_i E_i+K_Y-E_i|$ o\`u $E_i$ d\'esigne la fibre
exceptionnelle de $\sigma_i$ et $m_i\in\{0,1\}$ la multiplicit\'e de
$\wtilde{C}_i$ en $O_i$; par cons\'equent $\sigma_i^* D_i+(1-m_i)E_i\in
|\wtilde{C}+K_Y|$ et donc $D_i\in (\sigma_{i})_*|\wtilde{C}+K_{Y}|$. En
d\'ecomposant $\sigma_i$ en une succession d'\'eclatements, on d\'emontre ainsi que
\[(\sigma_i)_*|\wtilde{C}+K_Y|=|\wtilde{C}_i+K_{Y_i}|\] 
d'o\`u r\'esulte
\begin{eqnarray*}
(\pi_1)_*|\wtilde{C}_1+K_{Y_1}|&=&(\pi_1\sigma_1)_*|\wtilde{C}+K_Y|\\
&=&(\pi_2\sigma_2)_*|\wtilde{C}+K_Y|\\
&=&(\pi_2)_*|\wtilde{C}_2+K_{Y_2}|.
\end{eqnarray*}

Insistons sur le fait que si $h^0(\widetilde{C}+K_Y)\leq 1$, alors le
syst\`eme lin\'eaire adjoint de $C$ n'existe pas, par d\'efinition. En
g\'en\'eral, lorsqu'on parlera de syst\`eme lin\'eaire, il est sous-entendu
que celui-ci contient au moins un pinceau, \ie\ un syst\`eme lin\'eaire
de dimension 1.

\begin{pro}
Le syst\`eme adjoint $\adj(C)$ existe si et seulement si $g(C)>1$. De
plus, dans ce cas, l'intersection avec $C$ induit un isomorphisme
$\adj(C)\simeq
|K_{C}|$.
\label{pro1}
\end{pro}
Ici $g(C)$ d\'esigne le genre de la normalisation
$\wtilde{C}$ de $C$ et $K_C$ l'image par $\pi:\wtilde{C}\to C$ d'un
diviseur canonique sur $\wtilde{C}$.
\begin{proof}
Consid\'erons la suite exacte
\[\xymatrix{0\ar@{->}[r]&\calo_Y(-\widetilde{C})\ar@{->}[r]&\calo_Y\ar@{->}[r]&\calo_{\widetilde{C}}\ar@{->}[r]&0}\]

qui, apr\`es tensorisation avec $\calo_Y(\widetilde{C}+K_Y)$ induit la
suite exacte
\[\cdots\xymatrix{\cohz(Y,K_Y)\ar@{->}[r]&\cohz(Y,\widetilde{C}+K_Y)\ar@{->}[r]&\cohz(\widetilde{C},K_{\wtilde{C}})\ar@{->}[r]&\cohu(Y,K_Y)\ar@{->}[r]&\cdots}\]
(pour simplifier les notations, dans les groupes de cohomologie nous
identifions les diviseurs \`a leurs faisceaux associ\'es). Puisque $Y$
est rationnelle les deux termes extrêmes sont r\'eduits \`a $\{0\}$, d'o\`u
le r\'esultat.
\end{proof}

\begin{exe}
Soit $C\subset\plan$ une courbe irr\'eductible de degr\'e $d$, avec des
singularit\'es $p_1,\ldots,p_\ell$, de multiplicit\'es
$m_1,\ldots,m_\ell\geq 2$
respectivement ($\ell\geq 0$). Supposons que $g(C)\geq 2$; observons qu'alors $d$ est au moins \'egal \`a 4.

a) Si toutes les singularit\'es $p_1,\ldots,p_\ell$ sont ordinaires,
une r\'esolution plong\'ee des
singularit\'es de $C$ est obtenue en prenant pour $\pi:Y\to\plan$
l'\'eclatement
de $p_1,\ldots,p_\ell$ dans $\plan$. Alors $\adj(C)$ est le syst\`eme
lin\'eaire 
\[\pi_*|(d-3)L-\sum_{i=1}^\ell (m_i -1)E_i|\]
priv\'e de ses composantes fixes,
o\`u $L$ d\'esigne l'image inverse par $\pi$ d'une droite g\'en\'erale et
$E_i$ est la
courbe exceptionnelle au dessus de $p_i$. C'est le syst\`eme lin\'eaire
des
courbes de degr\'e $d-3$ qui passent par $p_i$ avec multiplicit\'e
$m_i-1$ ($i=1,\ldots,\ell$), priv\'e de ses composantes fixes. 

b) Supposons que $C$ poss\`ede des singularit\'es non ordinaires. Toute
r\'esolution plong\'ee $\pi$ de $C$ se factorise \`a travers l'\'eclatement de
$p_1,\ldots,p_\ell$ dans $\plan$, mais elle peut être diff\'erente de
celui-ci. N\'eanmoins, si la r\'esolution est minimale, $\adj(C)$ est un
syst\`eme
lin\'eaire de la forme
\[\pi_*\Big|(d-3)L-\sum_{i=1}^\ell (m_i -1)E_i-\sum a_{ij} F_{ij}\Big|\]
priv\'e de ses composantes fixes, o\`u $F_{ij}$ est un diviseur effectif contract\'e par $\pi$
sur $p_i$
et $a_{ij}> 0$. 

\label{exe1}
\end{exe}

\begin{rem}
Comme il suit de l'exemple pr\'ec\'edent, lorsque $C$ est une courbe plane, le degr\'e des \'el\'ements de
$\adj(C)$ est
$\leq \deg(C)-3$. 
\label{rem*}
\end{rem}

\begin{exe}
a) Soit $C=C_{g+2}\subset\plan$ une courbe g\'en\'erale de degr\'e $g+2\geq 2$
avec un point singulier ordinaire $p$ de multiplicit\'e $g$ (donc
$g(C)=g$). Si $g\leq 1$ le syst\`eme $\adj(C)$ n'existe pas et si
$g\geq 2$, il est constitu\'e de $(g-1)$ droites passant par $p$.

b) Si $C\subset\plan$ est une sextique avec deux points triples
ordinaires $p$
et $q$, alors, $\adj(C)$ est obtenu \`a partir du syst\`eme des cubiques
avec point double en $p$ et $q$ en suprimant la droite $pq$ qui est
une composante fixe: ainsi $\adj(C)$ est le syst\`eme lin\'eaire des
coniques passant par $p$ et $q$.

c) Le syst\`eme lin\'eaire $\cald$ des cubiques planes passant par un
ensemble $I$ de 7 points en position g\'en\'erale d\'efinit une application
rationnelle $\gamma:\proj^2\tor\plan$ de degr\'e 2, l'involution
correspondante $\sigma$ \'etant une involution de Geiser. Consid\'erons
l'\'eclatement $\pi:X\to\plan$ de $I$: alors
$\sigma':=\pi^{-1}\sigma\pi$ est un automorphisme de $X$ (c.f.
\cite{BayBea}). Puisque le groupe des classes de diviseurs sur $X$ invariants par $\sigma'$ est engendr\'e par $K_X$, la courbe $C$ des points fixes de
$\sigma$ est de degr\'e $3d$ avec multiplicit\'e $d$ en les points de
$I$. Par ailleurs, la restriction de $\gamma$ \`a un membre g\'en\'eral $D$
de $\cald$ est un revêtement double $D\to\proj^1$ ramifi\'e en 4 points
puisque $g(D)=1$. Ainsi l'intersection libre de $C$ et $D$ est
constitu\'ee de 4 points, d'o\`u aussit\^ot $d=2$: la courbe $C$ des points
fixes de l'involution de Geiser est donc une sextique avec points
doubles ordinaires en les 7 points de $I$ et $\adj(C)$ est le syst\`eme
lin\'eaire des cubiques passant par $I$. 

d) Le syst\`eme lin\'eaire $\cald$ des sextiques planes singuli\`eres sur
un ensemble $J$ de 8 points en position g\'en\'erale d\'efinit une
application rationnelle $\gamma:\plan\tor Q\subseteq \proj^3$ de degr\'e 2, o\`u $Q$ est un c\^one quadratique, l'involution correspondante \'etant
une involution de Bertini $\sigma$. Comme ci-dessus on observe que la
courbe des points fixes de $\sigma$ est de degr\'e $3d$ avec
multiplicit\'e $d$ sur $J$. La restriction de $\gamma$ \`a un membre
g\'en\'eral $D$ de $\cald$ est un revêtement double $D\to\proj^1$ ramifi\'e
en 6 points puisque $g(C)=2$ et par cons\'equent $d=3$: la courbe $C$
des points fixes de l'involution de Bertini est une nonique avec points triples en
les 8 points de $J$ et $\adj(C)$ est le syst\`eme lin\'eaire des
sextiques singuli\`eres sur $J$. 

\label{exe2}
\end{exe} 

Soit $\varphi:X_1\tor X_2$ une application birationnelle entre deux
surfaces (lisses et rationnelles). Si $D\subset X_1$
est une courbe dont aucune composante n'est contract\'ee par $\varphi$,
nous  notons
$\widetilde{\varphi}(D)$ la transform\'ee stricte directe de $D$ par
$\varphi$, et, si $\calh$ est
un syst\`eme lin\'eaire sur $X_1$ sans composante fixe,
$\widetilde{\varphi}(\calh)$
d\'esigne le syst\`eme lin\'eaire sur $X_2$ engendr\'e par les
$\widetilde{\varphi}(D)$ o\`u $D\in\calh$ est un \'el\'ement g\'en\'eral de
$\calh$. Par construction $\widetilde{\varphi}(\calh)$ n'a pas de
composantes fixes; il s'appelle le \emph{transform\'e homalo\"idal} de $\calh$
par $\varphi$. On peut aussi le d\'efinir via une r\'esolution de
l'ind\'etermination de $\varphi$
\begin{equation}
\xymatrix{ &Y\ar@{->}[dl]_{\sigma_1}\ar@{->}[dr]^{\sigma_2}\\
X_1\ar@{-->}[rr]^\varphi& &X_2}
\label{triangle}
\end{equation}
o\`u $\sigma_1, \sigma_2$ sont des morphismes birationnels:
$\widetilde{\varphi}(\calh)$ est le syst\`eme lin\'eaire
$(\sigma_2)_*\sigma_1^*\calh$ priv\'e de ses \'eventuelles composantes
fixes. 

La propri\'et\'e fondamentale du syst\`eme adjoint est sa ``covariance"
relativement aux transformations de
Cremona; elle joue un r\^ole important dans la litt\'erature classique
lors de l'\'etude de l'op\'eration du groupe de Cremona dans les syst\`emes
lin\'eaires de courbes planes (\cite[chap. VIII]{God}, \cite[Book IV,
chap. VII]{Coo},  \cite[Libro quinto Cap.II]{ECH}).

\begin{pro}
Soient $\varphi:X_1\tor X_2$ une application birationnelle entre deux
surfaces et $X_1\supseteq C_1$ une courbe irr\'eductible telle que
$g(C_1)>1$. Alors
\[\widetilde{\varphi}(\adj(C_1))=\adj(\widetilde{\varphi}(C_1)).\]
\label{pro2}
\end{pro}
\begin{proof}
Consid\'erons un triangle commutatif comme dans le diagramme
(\ref{triangle}); 
sans perdre de g\'en\'eralit\'e, on peut supposer que $\sigma_1$ r\'esoud les
singularit\'es de $C_1$ et donc que $\sigma_2$ r\'esoud celles de
$\widetilde{\varphi}(C_1)$. (On observe que de l'hypoth\`ese $g(C_1)>1$
suit que $\widetilde{\varphi}(C_1)$ est bien d\'efinie). Notons
$\wtilde{C}$ la transform\'ee stricte inverse de $C_1$ par $\sigma_1$
(ou de $\widetilde{\varphi}(C_1)$ par $\sigma_2$). Par d\'efinition, \`a
composantes fixes pr\`es, on a
\[(\sigma_1)_*|\widetilde{C}+K_Y|=\adj(C),\ \
(\sigma_2)_*|\widetilde{C}+K_Y|=\adj(\widetilde{\varphi}(C_1)),\]
d'o\`u suit l'assertion.
\end{proof}

Soit $\calh$ un syst\`eme lin\'eaire sur une surface $X$. Nous noterons
$\alpha_\calh:X\tor \calh^{\vee}$ l'application rationnelle associ\'ee
\`a $\calh$. Dans le cas o\`u $\calh=\adj(C)$ on \'ecrira $\alpha_C$ au
lieu de $\alpha_{\adj(C)}$.

Soit $X\supseteq C$ une courbe irr\'eductible avec $g(C)>1$; d'apr\`es la
proposition \ref{pro1} on a le diagramme commutatif
\[\xymatrix{X\ar@{-->}[rr]^{\alpha_C}& &\adj(C)^{\vee}\\
& C\ar@{^{(}->}[lu]\ar@{-->}[ru]_{\gamma_C}&}\]
o\`u $\gamma_C$ est ``l'application canonique" de $C$, \ie\ que
$\gamma_C$ pr\'ec\'ed\'ee de la normalisation $\wtilde{C}\to C$ est le
morphisme canonique de $\wtilde{C}$

Le r\'esultat suivant est une cons\'equence directe de la proposition
\ref{pro2}.

\begin{cor}
Soient $\varphi:X_1\tor X_2$ une application birationnelle entre deux
surfaces et $X_1\supseteq C_1$ une courbe irr\'eductible telle que
$g(C_1)>1$; notons $C_2=\wtilde{\varphi}(C_1)$. Alors, il existe un
isomorphisme $\phi:\adj(C_1)^{\vee}\to \adj(C_2)^{\vee}$  qui rend
commutatif le diagramme suivant
\[\xymatrix{&X_1\ar@{-->}[rr]^{\varphi}\ar@{-->}[dd]^{\alpha_{C_1}}&
& X_2\ar@{-->}[dd]^{\alpha_{C_2}}\\
C_1\ar@{^{(}->}[ru]\ar@/_/@{-->}[rr]\ar@{-->}[rd]^{\gamma_{C_1}}&
&C_2\ar@{^{(}->}[ru]\ar@{-->}[rd]^{\gamma_{C_2}}& \\
&\adj(C_1)^{\vee}\ar@{->}[rr]^\phi&  &\adj(C_2)^{\vee}.
}\]
\vspace{-0.8 cm}

\qed
\label{cor3}
\end{cor}


On dit qu'un syst\`eme lin\'eaire sur une surface $X$ est
\emph{r\'eductible} s'il est compos\'e de diviseurs r\'eductibles; sinon,
on dit qu'il est \emph{irr\'eductible}.

\begin{pro}
Soit $\calh$ un syst\`eme lin\'eaire irr\'eductible sur une surface $X$.
Alors $\adj(H)$ ne d\'epend pas du choix de $H$ g\'en\'eral dans $\calh$.
\label{pro33}
\end{pro}
\begin{proof}
En effet, apr\`es un \'eclatement convenable et une suppression
\'eventuelle de composantes fixes, on peut supposer que $\calh$ est
sans points-base; dans ce cas les membres g\'en\'eraux de $\calh$ sont
lisses d'apr\`es le th\'eor\`eme de Bertini \cite[Chap. III, Cor.
10.9]{Har}. 
\end{proof}

Cette proposition permet de d\'efinir le \emph{syst\`eme adjoint}
$\adj(\calh)$ d'un syst\`eme lin\'eaire irr\'eductible $\calh$: il s'agit
simplement de $\adj(H)$ pour $H$ g\'en\'eral dans $\calh$.

Le r\'esultat ci-dessous suit imm\'ediatement de la proposition
\ref{pro2}.

\begin{cor}
Soit $M$ un groupe de transformations birationnelles de la surface
$X$ et $\calh$ un syst\`eme lin\'eaire irr\'eductible sur $X$. Si $\calh$
est stable par $M$ (\ie\ $\wtilde{\mu}(H)\in\calh$ pour tout $H$ g\'en\'eral
dans $\calh$ et $\mu\in M$), alors $\adj(\calh)$ est aussi stable par
$M$.
\label{cor33}\qed
\end{cor}

Dans la terminologie classique, on dit qu'un syst\`eme lin\'eaire $\call$
sans composantes fixes et r\'eductible est \emph{compos\'e avec un
pinceau} $\call^c$: rappelons ce dont il s'agit.

\begin{lem}
Soit $\call$ un syst\`eme lin\'eaire sans composantes fixes et r\'eductible.
Alors il existe un unique pinceau irr\'eductible $\call^c$ tel que tout
\'el\'ement de $\call$ est compos\'e d'une somme d'\'el\'ements de $\call^c$;
l'application $\alpha_\call$ associ\'ee \`a $\call$ se factorise en 
\[\xymatrix{X\ar@{-->}[rr]^{\alpha_{\call^c}}\ar@{-->}[dr]_{\alpha_\call}&
& ({\call^c})^{\vee}\simeq\proj^1\ar@{->}[ld]\\
&\call^{\vee}&}.\]
\label{lem3}
\end{lem}
\begin{proof}
Quitte \`a remplacer $X$ par un \'eclat\'e convenable, on peut supposer que
$\alpha_\call$ est un morphisme. D'apr\`es le th\'eor\`eme de Bertini
(\cite[\S 7.9]{Iit}), $\alpha_\call(X)$ est une courbe.

Consid\'erons maintenant la factorisation de Stein de $\alpha_\call$
\[\xymatrix{X\ar@{->}[rr]^{\beta}\ar@{->}[dr]_{\alpha_\call}& &
W\ar@{->}[ld]\\
&\alpha_\call(X)&}\]
o\`u $W$ est une courbe normale (donc lisse) et les fibres de $\beta$
sont connexes (\cite[\S 2.13 et 2.14]{Iit}) et en g\'en\'eral lisses
(\cite[\S 7.9]{Iit}). Puisque $X$ est rationnelle, $W$ est isomorphe
\`a $\proj^1$ de sorte que les fibres de $\beta$ constituent le pinceau
cherch\'e.
\end{proof}

\begin{lem}
Soit $\calh$ un syst\`eme lin\'eaire irr\'eductible sur la surface $X$.
Alors tous les membres de $\calh$ sont connexes (\ie\ que leurs
supports le sont).
\label{lem33}
\end{lem}
\begin{proof}
Lorsque $\dim\,\alpha_{\calh}(X)=2$ le r\'esultat est bien connu: voir
\cite[Chap. I. \S 3.3]{Laz}.

Supposons donc que $\alpha_{\calh}(X)$ soit une courbe: par
irr\'eductibilit\'e, $\calh$ est un pinceau et alors en utilisant la
factorisation de Stein de $\alpha_{\calh}$, on voit que les fibres de
$\alpha_{\calh}$, \ie\ les \'el\'ements de $\calh$, sont connexes.
\end{proof}

Soit $\calh$ un syst\`eme lin\'eaire irr\'eductible sur $X$. On se propose
de d\'efinir les \emph{adjoints successifs} de $\calh$. Supposons donc
que $\adj(\calh)$ existe; s'il est irr\'eductible on pose
$\adj(\calh)^{(1)}=\adj(\calh)$ et sinon
$\adj(\calh)^{(1)}=\adj(\calh)^c$. Ensuite on d\'efinit
$\adj(\calh)^{(n)}=\adj\left(\adj(\calh)^{(n-1)}\right)^{(1)}$. Par
construction les syst\`emes lin\'eaires $\adj(\calh)^{(i)}$ sont
irr\'eductibles.

\begin{pro}
Soit $\calh$ un syst\`eme lin\'eaire irr\'eductible sur $X$. Il existe
$d\geq 0$ tel que $\adj(\calh)^{(i)}$ n'existe pas pour $i>d$.
\label{pro34}
\end{pro}

\begin{proof}
Soit $\varphi:X\tor\plan$ une application birationnelle. D'apr\`es la
proposition \ref{pro2}, en tenant compte du lemme \ref{lem3}, si
$\adj(\calh)^{(i)}$ existe alors
\[\wtilde{\varphi}(\adj(\calh)^{(i)})=\adj(\wtilde{\varphi}(\calh))^{(i)},\]
ce qui ram\`ene au cas o\`u $X=\plan$. Mais nous savons (remarque
\ref{rem*}) que $\deg\,\adj(\calh)^{(1)}\leq \deg\,\calh-3$ lorsque
$X=\plan$, d'o\`u aussit\^ot le r\'esultat.
\end{proof}

\begin{pro}
Soit $M$ un groupe de transformations birationnelles de $X$.
Supposons qu'il existe une courbe irr\'eductible $C$ de genre $>1$
telle que $\wtilde{\mu}(C)=C$ pour tout $\mu\in M$. Il existe alors
un syst\`eme lin\'eaire irr\'eductible $\calh$ sur $X$, stable par $M$ et
tel que $g(H)\leq 1$ pour tout $H$ g\'en\'eral dans $\calh$.
\label{pro35}
\end{pro}
\begin{proof}
Partant de $C$ on construit la suite des adjoints $\adj(C)^{(1)}$,
$\adj(C)^{(2)}$, $\ldots$ dont le dernier convient d'apr\`es la
proposition \ref{pro1}.
\end{proof}

\begin{exes}

a) Soit $C\subseteq X$ une courbe hyperelliptique de genre $>1$ constitu\'ee de points fixes d'un automorphisme $\sigma$ qui est une involution de de Jonqui\`eres. On sait que la paire $(X,C)$ est birationnellement \'equivalente \`a $(\plan,C_{g+2})$, o\`u $C_{g+2}$ est une courbe irr\'eductible de degr\'e $g+2$ avec un point singulier ordinaire $g$-uple (c.f. \cite{BayBea}). Par cons\'equent $\adj(C)^{(1)}$ est un pinceau de courbes rationnelles. 

Consid\'erons maintenant le groupe $M_C$ des transformations birationnelles $\mu$ de $X$ telles que $\wtilde{\mu}(C)=C$: ce groupe stabilise le pinceau $\adj(C)^{(1)}$ et par cons\'equent est conjugu\'e \`a un sous-groupe de $\jon(\plan)$.

On peut aussi raisonner directement: quitte \`a effectuer une contraction, on peut supposer que $(X,\sigma)$ est minimale; alors il existe un fibr\'e en coniques $r:X\to\proj^1$ tel que $r\circ\sigma=r$. De plus le sous-groupe des invariants de $\sigma$ dans le groupe de Picard de $X$ est de rang 2 engendr\'e par $K_X$ et la classe $f$ d'une fibre de $r$; on a donc $C\equiv aK_X+bf$ avec $a=-1$ puisque $C\cdot f=2$ et $K_X\cdot f=-2$; par cons\'equent $\adj(C)^{(1)}=|f|$ et le pinceau cherch\'e n'est rien d'autre que le pinceau associ\'e \`a $r$.

b) Consid\'erons le syst\`eme bi-anticanonique $|-2 K_S|$ d'une surface de Del Pezzo $S$ de degr\'e 1: un membre g\'en\'eral $D$ de ce syst\`eme est de genre 2 et donc est hyperelliptique; de plus $\adj(|-2 K_S|)=|-K_S|$ est un pinceau de courbes dont l'\'el\'ement g\'en\'eral est de genre 1.
\label{exes4}
\end{exes}

La th\'eorie des adjoints montre ici que la paire $(X,C)$ de (a), lorsque $g=2$,  n'est pas birationnellement \'equivalente \`a la paire $(S,D)$ de (b).

\begin{rem}
Partant d'un syst\`eme lin\'eaire sur $X$ qui est stable par $M$, on a la
même conclusion que dans la proposition \ref{pro35}. Dans la litt\'erature classique, ce r\'esultat est la cl\'e
pour la classification des sous-groupes ``continus" du groupe de
Cremona: voir \cite[Book IV, chap. VIII]{Coo}, \cite[Libro quinto, \S
22]{ECH}. Mais cela s'applique aussi aux groupes finis: si $M$ est
une sous-groupe fini d'automorphismes de $X$, l'image r\'eciproque par
l'application canonique $X\to X/M$ d'un syst\`eme lin\'eaire ample, par
exemple, est irr\'eductible et stable par $M$. Avec la m\'ethode des adjoints successifs,
on montre que $M$ stabilise un syst\`eme lin\'eaire irr\'eductible dont les
\'el\'ements g\'en\'eraux sont de genre $\leq 1$; cette r\'eduction est
semblable \`a celle \`a laquelle on a fait allusion dans l'introduction
et qui est issue de la th\'eorie des contractions des rayons extremaux
du c\^one $\overline{NE(X)}$ des courbes sur $X$ (c.f. \cite[\S
2.18]{CCKM}).
\end{rem}

Soit $\varphi:X\tor X$ une application birationnelle. Un point $x\in
X$ est un \emph{point fixe} de $\varphi$ si $\varphi$ est d\'efinie en
$x$ et $\varphi(x)=x$: on notera $\fix(\varphi)$ l'adh\'erence dans $X$
de l'ensemble des points fixes par $\varphi$.

On dit qu'une application rationnelle $p:X\tor\proj^1$ est une
\emph{fibration rationnelle} (resp. \emph{elliptique}) si la fibre
g\'en\'erale de $p$ est une courbe rationnelle (resp. elliptique)

Le r\'esultat suivant constitue la premi\`ere \'etape de la preuve du
th\'eor\`eme de Castelnuovo.

\begin{cor}
Soit $M$ un groupe de transformations birationnelles de la surface
$X$. Supposons qu'il existe une courbe irr\'eductible $C$ de genre $>1$
telle que $C\subset \fix(\mu)$ pour tout $\mu\in M$. Il existe alors
une fibration rationnelle ou elliptique $p:X\tor\proj^1$ telle que
$p\circ\mu=p$ pour tout $\mu\in M$. 
\label{cor34}
\end{cor}

\begin{proof}
Consid\'erons un syst\`eme lin\'eaire $\calh$ comme dans la proposition et
l'application rationnelle $\alpha_{\calh}:X\tor\calh^{\vee}$
correspondante qui est \'equivariante. En utilisant le lemme ci-dessous,
on observe que $\alpha_{\calh}(C)$ engendre $\calh^{\vee}$ puisque
aucun \'el\'ement de $\calh$ ne peut contenir $C$ qui est de genre $>1$
par hypoth\`ese. Maintenant, comme $C$ est fixe, le groupe $M$ op\`ere
trivialement au but, d'o\`u l'assertion.
\end{proof}

\begin{lem}
Soient $\calh$ un syst\`eme lin\'eaire irr\'eductible sur $X$ et $C$ une
courbe irr\'eductible contenue dans le support d'un \'el\'ement de $\calh$.
Alors le genre d'un \'el\'ement g\'en\'eral de $\calh$ est $\geq g(C)$. 
\label{lem34}
\end{lem}
\begin{proof}
Quitte \`a effectuer des \'eclatements, on peut supposer que $\calh$ est
sans points-base et que $C$ est lisse; en particulier un \'el\'ement
g\'en\'eral $H$ de $\calh$ est aussi lisse. De la suite exacte
\[\xymatrix{0\ar@{->}[r]&\calo(C+K_X)\ar@{->}[r]&\calo(H+K_X)}\]
il suit que $\dim\,|C+K_X|\leq \dim\,|H+K_X|$; de plus, la
restriction de ces deux syst\`emes lin\'eaires \`a $C$ et $H$
respectivement, d\'efinissent les s\'eries canoniques de ces deux
courbes, d'o\`u l'assertion.
\end{proof} 

\begin{exe}
Consid\'erons l'involution de Bertini (exemple \ref{exe2}d) dont la
courbe $C$ des points fixes est une nonique avec points triples
situ\'es sur un ensemble $J$ de 8 points en position g\'en\'erale. Ici le
pinceau cherch\'e est $\adj(C)^{(2)}$ qui est constitu\'e des
cubiques passant par $J$.
\end{exe}

\section{Preuve du th\'eor\`eme de Castelnuovo}
\label{sec:PreuveCast}
Soit $F$ une transformation birationnelle de $\plan$ qui n'est pas
l'identit\'e. Supposons que $\fix(F)$ contienne une courbe $C$ de genre
$>1$. On sait qu'il existe une fibration rationnelle ou elliptique
$p:\proj^2\tor\proj^1$ telle que $F\circ p=p$ (corollaire
\ref{cor34}). Soit $\sigma:X\to\proj^2$ un morphisme birationnel tel
que $q:=p\circ\sigma$ est un morphisme; posons
$G:=\sigma^{-1}F\sigma$. 

Si $p$ est rationnelle, d'apr\`es le th\'eor\`eme de Noether-Enriques
(\cite[Thm. III.4]{Bea}) il existe un ouvert $U\subseteq\proj^1$ tel
que 
\[\xymatrix{q^{-1}(U)\ar@{->}[rr]^{\sim}\ar@{->}[dr]_{q}& &
U\times\proj^1\ar@{->}[ld]^{pr_1}\\
&U&}\]
d'o\`u suit que $F$ est conjugu\'ee \`a une transformation de de
Jonqui\`eres. Comme la courbe $C$ n'est pas rationnelle et $F\neq id$,
on en d\'eduit que la restriction de $G$ \`a l'une des fibres g\'en\'erales
de $q$ fixe exactement deux points: il s'ensuit que $q$ pr\'esente $C$
comme revêtement \`a deux feuilles de $\proj^1$, et donc que $C$ est
hyperelliptique. De plus, il existe $a\in{\rm PGL}(2,\complex(x))$
tel que $F$ s'exprime dans des coordonn\'ees affines convenables par
\[F_a: (x,y)\mapstor
\left(x,\frac{a_{11}(x)y+a_{12}(x)}{a_{21}(x)y+a_{22}(x)}\right);\]
ici $A=(a_{ij}(x))$ repr\'esente $a$.

Si $(y^2=h(x))$ est l'\'equation affine de $C$, on v\'erifie directement
que $F_a$ fixe $C$ si et seulement si $A=\left(\begin{array}{ll}
a_1&ha_2\\a_2&a_1\end{array}\right)$. 

Le lemme suivant ach\`eve de d\'emontrer le th\'eor\`eme lorsque la fibration
est rationnelle. 

\begin{lem}
Si $a$ (ou $F_a$) est d'ordre fini et si $F_a$ fixe une courbe de
genre $>0$, alors $a$ (ou $F_a$) est une involution.
\label{leminv}
\end{lem}
\begin{proof}
On observe que la condition sur l'ensemble des points fixes de $F_a$
implique que $A$ n'est pas diagonalisable et donc que le polyn\^ome
minimal $m_A$ de $A$ est irr\'eductible dans $\C(x)[T]$.

Notons $n$ l'ordre de $a$: il existe donc $g\in\complex(x)$ tel que
$A^n=gI$ ($I$ d\'esigne la matrice identit\'e) et par suite on a
$(\det\,A)^n=g^2$.

Si $n=2d$ est pair, alors $g=f^d$ pour un \'el\'ement $f\in\complex(x)$.
Maintenant $m_A$ divise $T^{2d}-g=T^{2d}-f^d$: on a donc $m_A=T^2-\xi
f$ o\`u $\xi$ est une racine $d^{\mbox{i\`eme}}$ de l'unit\'e, c'est-\`a-dire
que $a$ est d'ordre 2.

Si $n$ est impair, $g=f^n$ pour un \'el\'ement $f\in\complex(x)$. En
effet, choisissons deux entiers $u$ et $v$ tels que $nu+2v=1$; on a
alors
\[g=g^{un}(g^2)^v=(g^u(\det\,A)^v)^n.\]
Cette fois $T^n-g=T^n-f^n$ se d\'ecompose dans $\complex(x)[T]$ ce qui
n'est pas possible puisque $m_A$ est irr\'eductible et divise $T^n-g$.
\end{proof}

Supposons maintenant que $p$, et donc $q$, est elliptique. La
restriction $G_b$ de $G$ \`a une fibre g\'en\'erale $X_b$ de $q$ est un automorphisme
admettant au moins $2$ points fixes: en effet, la courbe de points
fixes $C$ \'etant de genre $>1$ n'est pas contenue dans une fibre de
$q$ et la restriction de $q$ \`a $C$ n'est pas birationnelle. Dans ces
conditions, on sait que $G_b$ est d'ordre fini \'egal \`a $2,3,4$ ou $6$:
voir \cite[Chap. IV, \S 4.7]{Har}; il s'ensuit que $G$ lui-même est
d'ordre $2,3,4$ ou 6, ce qui permet de supposer, si c'est plus
confortable, que $G$ est un automophisme de $X$.

L'existence de $2$ points fixes dans $X_b$ exclut l'ordre 6: d'apr\`es
\cite[Chap. IV, \S 4.20.2]{Har}, si $G$ est d'ordre 6, alors $X_b$
est isomorphe au quotient de $\C$ par $\Z[\omega]$ o\`u $\omega=e^{2\pi
i/3}$ et $G_b$ est la multiplication par $-\omega$ ou $-\omega^2$ qui
n'a qu'un seul point fixe. On observe par contre que l'ordre 3 n'est
pas exclu.

Si l'ordre de $G$ est 4, alors la fibre g\'en\'erale $X_b$ est isomorphe
au quotient de $\C$ par $\Z[i]$, o\`u $i^2=-1$, et $G_b$ est la multiplication par $\pm i$ qui poss\`ede exactement 2 points
fixes: voir \cite[Chap. IV, \S 4.20.1]{Har}. Ainsi la courbe de points fixes $C$ est hyperelliptique. On observe que $C$ est aussi une courbe de points fixes de l'involution $G^2$, qui est donc de de Jonqui\`eres, de sorte qu'on se trouve dans la situation \'etudi\'ee dans l'exemple \ref{exes4}a): il existe une fibration rationnelle $r$ qui est laiss\'ee stable par $G$. Puisque $g(C)>1$, la restriction de $r$ \`a $C$ est surjective et alors $G\circ r=r$ ce qui ram\`ene au cas \'etudi\'e plus haut; ceci montre que $G$ est d'ordre 2 .  Cette contradiction implique que l'ordre 4 n'est pas possible.

Consid\'erons maintenant un groupe $M$ de transformations de $\plan$
dont tous les \'el\'ements fixent une même courbe de genre $>1$, comme
dans le th\'eor\`eme. On sait alors que $M$ fixe une fibration rationnelle ou elliptique. Dans le premier cas $M$ est conjugu\'e \`a un sous-groupe du groupe de de Jonqui\`eres. Dans le second cas, puisque en g\'en\'eral le groupe des automorphismes d'une courbe elliptique fixant un point est cyclique d'ordre 2,4 ou 6, le groupe $M$ lui-même est cyclique; l'analyse pr\'ec\'edente montre alors que seuls les ordres 2 et 3 sont possibles. 
Les pr\'ecisions concernant la classe de
conjugaison de $M$ dans le groupe de Cremona de $\plan$ suivent alors
de \cite{BayBea}, \cite{Fer} ou \cite{Bla}.

\section{Remarques finales}
Dans cette derni\`ere section, nous d\'emontrerons le r\'esultat suivant, qui prouve que les th\'eor\`emes \ref{teo1} et \ref{teo2} ne se g\'en\'eralisent pas aux courbes rationnelles ou elliptiques (de genre $0$ ou $1$).
\begin{pro}
Soit $C$ une courbe plane qui est l'image d'une droite ou d'une cubique lisse par une transformation birationnelle de $\plan$. Le groupe des transformations de Cremona qui fixent points par points la courbe $C$ n'est ni d'ordre fini, ni ab\'elien, ni conjugu\'e \`a un sous-groupe de $\jon(\plan)$.
\label{propellrat}
\end{pro}
\begin{rem}
Ce groupe a \'et\'e introduit dans \cite{Giz} sous le nom de "groupe d'inertie de $C$". 
\end{rem}
\begin{rem}
L'hypoth\`ese de l'\'enonc\'e impose une restriction sur la courbe: en effet, il existe des courbes rationnelles planes qui ne sont pas l'image d'une droite par une transformation de Cremona, par exemple une sextique avec $10$ points doubles (voir \cite{Coo}, Book IV, chap. II, \S 2 ou \cite{KuMu}); de m\^eme, une sextique avec $9$ points doubles est une courbe de genre $1$ qui n'est pas l'image d'une cubique lisse par une transformation de Cremona.
Ceci peut \^etre v\'erifi\'e en observant que tout pinceau de courbes rationnelles planes intersecte une sextique avec uniquement des points doubles ordinaires en au minimum $4$ points en dehors des points-base (voir le lemme ci-dessous).
\end{rem}

\begin{lem}
Soit $C$ une sextique plane qui n'a que des doubles ordinaires. Soit $\Lambda$ un pinceau de courbes rationnelles planes.

Alors, $C$ intersecte une courbe g\'en\'erale de $\Lambda$ en au minimum $4$ points en dehors des points-base.
\end{lem}
\begin{proof}
Notons $n$ le degr\'e des courbes de $\Lambda$ et $m_1,...,m_k$ les multiplicit\'es des points-base (qui peuvent \^etre sur $\plan$ ou infiniment proches). La condition de rationalit\'e et le fait que le syst\`eme soit un pinceau donnent respectivement
\begin{eqnarray}\label{eqs1}\frac{(n-1)(n-2)}{2}-\sum_{i=1}^k \frac{m_i(m_i-1)}{2}&=&0\\
\label{eqs2}\frac{(n+1)(n+2)}{2}-\sum_{i=1}^k \frac{m_i(m_i+1)}{2}&=&2.\end{eqnarray}
En soutrayant l'\'equation (\ref{eqs1}) \`a  l'\'equation (\ref{eqs2}), on obtient \begin{equation}\label{eqs3}3n-\sum_{i=1}^k m_i=2.\end{equation}. 

En notant $P_1,...,P_l$ les points doubles de $C$ et respectivement $n_1,...,n_l$ les multiplicit\'es du pinceau $\Lambda$ en ces points (qui peuvent \^etre $0$ si les points ne sont pas des points-base), la courbe $C$ intersecte une courbe g\'en\'erale de $\Lambda$ en $6n-2\sum_{i=1}^l n_i=2(3n-\sum_{i=1}^l n_i)$ points en dehors des points-base. Comme $\sum_{i=1}^l n_i\leq \sum_{i=1}^k m_i$, le r\'esultat suit de l'\'equation (\ref{eqs3}).
\end{proof}

\bigskip

Citons quelques exemples: 

\begin{exe}
a) Soit $G$ le groupe des transformations lin\'eaires de $\plan$ du type \[(x:y:z) \mapsto (a x:y+bx:z+cx),\]
o\`u $a,b,c\in\C, a\not=0$. On voit directement que $G$ fixe (points par points) la droite $L$ de $\plan$ d'\'equation $x=0$ (en fait, tout automorphisme lin\'eaire de $\plan$ qui fixe $L$ appartient \`a $G$) et que $G$ n'est ni fini ni ab\'elien.
Comme $G$ est un groupe de transformations lin\'eaires et comme ses seuls points fixes sont sur la droite $L$, les pinceaux de droites invariants par $G$ sont constitu\'es des droites passant par un point de $L$. En d'autres termes, $G$ est un sous-groupe de $\jon_p(\plan)$ si et seulement si $p\in L$.

b) Notons de plus $H$ le groupe des transformations birationnelles du plan (exprim\'ees en coordonn\'ees affines dans l'ouvert $z=1$) du type
\[(x,y)\dasharrow \Big(\frac{x}{\alpha(y)x+\beta(y)},y\Big),\]
o\`u $\alpha(y),\beta(y) \in \complex(y)$ sont des fonctions rationnelles et $\beta\not=0$. On observe \`a nouveau que $H$ fixe la droite $L$ (d'\'equation affine $x=0$) et que $H$ n'est ni fini ni ab\'elien. On voit que $H$ laisse invariant le pinceau des droites de $\plan$ passant par le point $(1:0:0)$ (d'\'equations affines $y= a$, $a\in \C$), mais qu'aucun des autres pinceaux de droites du plan n'est invariant par $H$. En d'autres termes, $H$ est un sous-groupe de $\jon_p(\plan)$ si et seulement si $p$ est le point $(1:0:0)\notin L$.

c) Il suit des observations pr\'ec\'edentes que le groupe engendr\'e par $G$ et $H$ fixe points par points la droite $L$; qu'il n'est ni fini, ni ab\'elien et qu'il n'est pas un sous-groupe de $\jon_p(\plan)$, quel que soit le point $p\in\plan$. On peut \'egalement observer que l'intersection des groupes $G$ et $H$ est triviale. 
\label{exedroite}
\end{exe}

\begin{exe}
Soit $C$ une courbe cubique plane lisse. Pour tout point $p$ de $C$; notons $\sigma_p$ l'involution de centre $p$ qui fixe $C$ d\'efinie comme suit: si $D$ est une droite g\'en\'erale passant par $p$, on a $\sigma_p(D)=D$ et la restriction de $\sigma_p$ \`a $D$ est l'involution de points fixes $(D\cap C)\backslash\{p\}$. 

Si $p$ et $p'$ sont deux points distincts en position g\'en\'erale de $C$, remarquons que le groupe engendr\'e par $\sigma_p$ et $\sigma_{p'}$ est un groupe infini non ab\'elien.
En effet, les restrictions de $\sigma_p$ et $\sigma_{p'}$ \`a la droite passant par $p$ et $p'$ sont deux involutions avec un unique point fixe en commun et leur restriction \`a cette droite engendre un tel groupe.
\label{execubique}
\end{exe}

\vspace{0.3 cm}

Prouvons maintenant, gr\^ace \`a ces deux exemples, la proposition \ref{propellrat}. \`A l'aide d'une transformation de Cremona, on se ram\`ene au cas o\`u $C$ est une droite ou une cubique lisse. Notons $M_C$ le groupe des transformations de Cremona qui fixent points par points la courbe $C$. Les exemples \ref{exedroite} et \ref{execubique} montrent que $M_C$ n'est ni fini ni ab\'elien. 

Supposons maintenant que $M_C$ est conjugu\'e \`a un sous-groupe de $\jon(\plan)$, ce qui implique que $M_C$ laisse invariant un pinceau de courbes rationnelles $\Lambda$. On note $\base(\Lambda)$ l'ensemble des points-base de ce pinceau.
De m\^eme qu'\`a la section \ref{sec:PreuveCast}, on observe que  $C$ est rationnelle: sinon elle coupe une courbe g\'en\'erale de $\Lambda$ en au moins $2$ points en dehors de $\base(\Lambda)$; comme le groupe des automorphismes de $\mathbb{P}^1$ qui fixent $2$ points est ab\'elien, cette situation ne se pr\'esente pas. On suppose alors que $C$ est la droite  
d'\'equation $x=0$.

On observe que le groupe $G$ des transformations lin\'eaires qui fixent $C$ (voir exemple \ref{exedroite}) op\`ere dans $\base(\Lambda)$; puisque les points de $C$ sont les seules orbites finies de $G$, les points-base de $\Lambda$ appartiennent \`a $C$ ou lui sont infiniment proches.

Pour tous $\mu,\nu \in \C$, non tous deux nuls,
l'application rationnelle 
\[\varphi_{\mu,\nu}:(x:y:z)\dasharrow (-x(\mu y+\nu z):y(x+\mu y+\nu z):z(x+\mu y+\nu z))\]
est une involution quadratique qui fixe points par points la droite $C$. Son syst\`eme lin\'eaire associ\'e $\Phi$
est l'ensemble des coniques de $\plan$ passant par les points $(1:0:0)$ et  $(0:\nu:-\mu)$ et qui sont tangentes en ce dernier point \`a la droite d'\'equation $x+\mu y+\nu z=0$.

En choissisant $\mu$ et $\nu$ de telle sorte que le point $(0:\nu:-\mu)$ ne soit pas un point-base de $\Lambda$, l'intersection de $\Lambda$ et $\Phi$ en dehors des points-base est $2n$, o\`u $n$ est le degr\'e des courbes de $\Lambda$. L'involution $\varphi_{\mu,\nu}$ envoie donc les courbes de $\Lambda$ sur des courbes de degr\'e $2n$, ce qui prouve que $\Lambda$ n'est pas invariant par $\varphi_{\mu,\nu}$ et donc que $M_C$ n'est pas birationnellement conjugu\'e \`a un sous-groupe de 
$\jon(\plan)$.

 \end{document}